\documentclass[12pt]{scrartcl}
\usepackage{latexsym}
\usepackage{amsmath,amssymb}
\usepackage{mathrsfs}
\usepackage{enumerate}
\usepackage{exscale}
\usepackage{parskip}

\newcounter{satzzaehler}[section]
\newcommand{\aussage}[1]
{\refstepcounter{satzzaehler}{\bf\thesection.\thesatzzaehler \ #1} \ }%

\def\scirc{\mathbin{\raise.15ex\hbox{\scriptsize$\circ$}}}
\def\eps{\varepsilon}
\def\qed{\hspace*{\fill} $\Box$\par\bigskip}

\title{Criteria for strong and weak random attractors}

\author{Hans Crauel%
  \thanks{Fachbereich Mathematik, 
          Johann Wolfgang Goethe-Universit\"at Frankfurt, 
          Robert-Mayer-Stra\ss e 10, 
          60325 Frankfurt, FRG; \ 
        \small{\tt crauel{\scriptsize @}math.uni-frankfurt.de}}
\and Georgi Dimitroff%
  \thanks{Institut Techno- und Wirtschaftsmathematik, 
  Fraunhofer-Platz~1, 67663 Kaiserslautern, FRG; \ 
  \small{\tt dimitroff{\scriptsize @}gmail.com}}
\and Michael Scheutzow%
  \thanks{Institut f\"ur Mathematik, MA 7-5, Fakult\"at II, 
        Technische Universit\"at Berlin, 
        Stra\ss e des 17. Juni 136, 10623 Berlin, FRG;  \ 
        \small{\tt ms{\scriptsize @}math.tu-berlin.de}}}
\date{}

\begin{document}  \maketitle

\begin{abstract}\noindent
  The theory of random attractors has different notions of 
  attraction, amongst them pullback attraction and weak 
  attraction. 
  We investigate necessary and sufficient conditions for the 
  existence of pullback attractors as well as of weak attractors. 
  \par\medskip

  \ Secondary\, 
37B25 \ 
37H99 \ 
37L55 \ 
60D05 \ 
\end{abstract}


\section{Introduction}\label{intro}
The notion of an attractor is one of the basic concepts in the 
theory of dynamical systems. 
For deterministic systems this notion has been of importance 
since several decades. 
Since about fifteen years also attractors for stochastic systems 
have been taken under consideration. 
The crucial obstacle came from the fact that the classical approach, 
using the Markov property of individual solutions to define the Markov 
semigroup and an associated generator, could not deal with joint motions 
of two or more points. 
This began to change after the introduction of stochastic flows, going 
back to Kunita and to Elworthy, and the introduction of the notion of 
random dynamical systems (see Arnold~\cite{arnold}). 
All approaches to random attractors use the theory of random dynamical 
systems (RDS). 
The notion of pullback attractors, which are referred to as strong 
attractors here, go back to Crauel, Flandoli, and 
Debussche~\cite{crauel-flandoli, crauel-debussche-flandoli}, 
and to Schmalfu\ss~\cite{schmalfuss94}. 
Later the notion of weak attractors was introduced by Ochs~\cite{ochs}. 
A comparison of these two concepts and yet another one, called 
forward attractor, has been investigated by 
Scheutzow~\cite{scheutzow02}, see also Crauel~\cite{crauel02a}. 
Several investigations, among others of notions of local random 
attractors, have been carried out by Ashwin and 
Ochs~\cite{ashwin-ochs}.  \par

The present paper gives necessary and sufficient criteria for 
the existence of strong and weak attractors for compact and 
for bounded sets. 
This is developed in Theorems~\ref{sq3}.\ref{strongB} 
and~\ref{sq3}.\ref{strongC} for strong attraction of bounded 
sets and compact sets, respectively, and in 
Theorems~\ref{weak}.\ref{weakB} and~\ref{weak}.\ref{weakC} 
for weak attraction of bounded sets and compact sets, respectively. 
A version of Theorem~\ref{weak}.\ref{weakB} for RDS on the 
Euclidean space~$\mathbb R^n$ goes back to Dimitroff~\cite{dimitroff}. 
Furthermore, it is shown that a weak attractor is a strong 
attractor if and only if it contains the Omega limit set of 
every set from the `domain of attraction'.  \par

Existence of strong attractors has been verified for many 
concrete systems in the literature. 
Also for the more general non-autonomous case existence of 
attractors has been obtained for a large variety of systems. 
There the `almost surely' used in the context of random 
attractors is replaced by `for all', while the method of the 
construction remains the same. 
See~\cite{crauel-debussche-flandoli} for this approach. 
The conditions for the existence of strong random attractors 
obtained in Theorems~\ref{sq3}.\ref{strongB} 
and~\ref{sq3}.\ref{strongC} are purely probabilistic. 
They have no non-autonomous version.  \par

One may wonder whether these conditions are of use in order to 
obtain existence of random attractors. 
In this respect we refer to recent work of Dimitroff and 
Scheutzow~\cite{dimitroff-scheutzow}, where the existence of 
strong random attractors for certain concrete systems is verified 
by using the conditions Theorem~\ref{sq3}.\ref{strongB} and 
Theorem~\ref{sq3}.\ref{strongC}.  \par

It should also be noted that a strong (or `pullback', see 
Remark~\ref{sq2}.\ref{tq7}) attractor allows for several 
explicit representations. 
It may be constructed by taking the union of all Omega limit 
sets of deterministic compact sets. 
Alternatively, the attractor is the Omega limit set of every 
deterministic compact set which contains the attractor with 
positive probability. 
Explicit constructions of weak attractors, on the other hand, 
have not been available yet. 
Theorems~\ref{weak}.\ref{weakB} and~\ref{weak}.\ref{weakC} 
give constructions of weak random attractors. 
The constructions are technically considerably more involved 
than those for strong attractors.

\section{Set-up and preliminaries}\label{sq2}
In this section we first introduce some basic notions of random dynamical 
systems (RDS), referring to Arnold~\cite{arnold} for a comprehensive 
presentation. 
Then we give a brief introduction to the concepts of strong and weak 
random attractors.  \par

Let~$E$ be a Polish space (i.e.\ a separable topological space which is 
metrisable by a complete metric), and let~$\mathscr B$ be its Borel 
$\sigma$-algebra. 
We will often assume~$E$ to be equipped with a metric inducing the 
topology, which will then be denoted by~$d$. 
Furthermore, $d$ will be assumed to be complete whenever the argument 
needs it. 
We have not tried to single out when this is the case. 
For $x\in E$ and for a subset~$A$ of $E$ we write 
$d(x,A)=\inf_{a\in A}d(x,a)$.  \par\smallskip

\aussage{Definition}
\begin{enumerate}[(a)]
\item 
  $\bigl(\Omega,\mathscr F,P,(\vartheta_t)_{t\in\mathbb R}\bigr)$ 
  is called a \emph{metric dynamical system (MDS)}, if 
  $(\Omega,\mathscr F,P)$ is a probability space, and the family of 
  mappings 
  $\bigl\{\vartheta_t:\Omega\to\Omega:t\in\mathbb R\bigr\}$ 
  satisfies 
  \begin{enumerate}[(i)]
  \item the mapping $(\omega,t)\mapsto\vartheta_t(\omega)$ 
    is 
    $(\mathscr F\otimes\mathscr B(\mathbb R),\mathscr F)$-measurable, 
  \item $\vartheta_{s+t}=\vartheta_s\scirc\vartheta_t$ for 
    every $s,t\in\mathbb R$, and $\vartheta_0=\mbox{Id}_\Omega$, 
    and 
  \item for each $t\in\mathbb R$, $\vartheta_t$ preserves the 
    measure~$P$. 
\end{enumerate}
\item A \emph{random dynamical system (RDS)} on the measurable space 
  $(E,\mathscr B)$ over the MDS 
  $\bigl(\Omega,\mathscr F,P,(\vartheta_t)\bigr)$ with time~$\mathbb R^+$ 
  is a mapping 
  \begin{displaymath}
    \varphi:[0,\infty)\times E\times\Omega \to E, \qquad
    (t,x,\omega)\mapsto \varphi(t,x,\omega)
  \end{displaymath}
  with the following properties: 
  \begin{enumerate}[(i)]
  \item $\varphi$ is 
    $\left(\mathscr B([0,\infty))\otimes
      \mathscr B\otimes\mathscr F,\mathscr B\right)$-measurable. 
  \item For all $s,t\in[0,\infty)$
    \begin{displaymath}
      \varphi(t+s,\omega)
      =\varphi(t,\vartheta_s\omega)\scirc\varphi(s,\omega)\qquad
        \mbox{for all}\ \omega\in\Omega, 
   \end{displaymath}
   and $\varphi(0,\omega)=\mbox{Id}_E$ for all $\omega\in\Omega$. 
 \end{enumerate}
 The RDS~$\varphi$ is called \emph{continuous} if, in addition, 
 \begin{itemize}
 \item[(iii)]the mapping $x\mapsto\varphi(t,x,\omega)$ is continuous 
   for all $(t,\omega)\in[0,\infty)\times\Omega$. 
 \end{itemize}
\end{enumerate}
The present definition of an RDS is formulated for continuous 
time only, i.e.\ $t\in\mathbb R$ for the base flow, and 
$t\in[0,\infty)$ for the RDS itself. 
The general definition also allows for discrete time.  \par

A map $B:\Omega\to 2^E$, where $2^E$ denotes the power set of~$E$, is 
said to be a \emph{random set} if its graph 
$\{(x,\omega):x\in B(\omega)\}\subset E\times\Omega$ is an element 
of the product $\sigma$-algebra $\mathscr B\otimes\mathscr F$. 
We will assume without further mentioning that a random set is 
nonempty almost surely. 
A map~$B:\Omega\to 2^E$ is called a \emph{closed random set} or 
a \emph{compact random set}, respectively, if 
$\omega\mapsto B(\omega)$~takes values in the closed subsets or 
in the compact subsets, respectively, of~$E$ and if, in addition, 
$\omega\mapsto d\bigl(x,B(\omega)\bigr)$ is measurable for every 
$x\in E$. 
A closed random set is always a random set, whereas a random set 
taking values in the closed subsets is a closed random set 
if $(\Omega,\mathscr F)$ is universally measurable. 
See, for instance, Hu and Papageorgiou~\cite{hu-papageorgiou} 
Chapter~2.1--2, also for further characterizations. 
In particular, a non-empty compact random set can also be 
characterized by considering it as a random variable taking values 
in the metric space of all compact subsets of~$E$, equipped with 
the Hausdorff metric, see below. 
%
\par

The following definition is (essentially) due to Crauel and 
Flandoli~\cite{crauel-flandoli}, see also Crauel, Debussche and 
Flandoli~\cite{crauel-debussche-flandoli} and 
Schmalfu\ss~\cite{schmalfuss92,schmalfuss94}.  \par\smallskip

\aussage{Definition}  \label{f_ibf_watr}
Let~$\varphi$ be an RDS on~$E$ over the MDS 
$\bigl(\Omega,\mathscr F,(\vartheta_t)_{t\in\mathbb R},P\bigr)$. 
Let $\mathcal B\subset2^E$ be an arbitrary subset of the 
power set of~$E$. 
A random set $\omega\mapsto A(\omega)$ is called a 
\emph{strong $\mathcal B$-attractor} for~$\varphi$ if 
\begin{enumerate}[(a)]
\item $\omega\mapsto A(\omega)$ is a compact random set. 
\item $A$ is \emph{strictly $\varphi$-invariant}, that is, for each 
  $t\geq0$ there exists a set $\Omega_t$ of full measure, such 
  that 
  $\varphi(t,\omega)(A(\omega))=A(\vartheta_t\omega)$ for all 
  $\omega\in\Omega_t$. 
\item $\displaystyle\lim_{t\to\infty}\sup_{x\in B}
  d\bigl(\varphi(t,\vartheta_{-t}\omega)(x),A(\omega)\bigr)=0$  
  almost surely for every $B\in\mathcal B$. 
\end{enumerate}
In particular, a strong $\mathcal B$-attractor is called 
\begin{itemize}
\item \emph{$B$-attractor} in case that~$\mathcal B$ is the 
  set of all bounded subsets of~$E$
\item \emph{$C$-attractor} in case that~$\mathcal B$ is the 
  set of all compact subsets of~$E$. 
\end{itemize}
It should be mentioned that the set of all bounded subsets of~$E$ 
depends on the choice of the metric~$d$ on~$E$, whereas this is 
not the case for the set of all compact subsets of~$E$. 

\aussage{Remark}
The notion of a~$\mathcal B$-attractor may be modified to allow 
also for random sets $\omega\mapsto B(\omega)$ by demanding 
instead of~(c) 
\begin{displaymath}
  \lim_{t\to\infty}\sup_{x\in B(\vartheta_{-t}\omega)}
  d\bigl(\varphi(t,\vartheta_{-t}\omega)(x),A(\omega)\bigr)=0 
  \qquad\mbox{almost surely for every}\ B\in\mathcal B. 
\end{displaymath}
This is of interest, in particular, when dealing with \emph{local} 
random attractors. 
For global attractors it does not yield much more generality, 
since whenever a general $\mathcal B$ comprises the set of all 
deterministic compact sets, then in case a $\mathcal B$-attractor 
exists also a $C$-attractor exists, it is unique, and 
the $\mathcal B$-attractor coincides with the $C$-attractor, 
see Crauel~\cite{crauel99}.  \par\medskip

It has turned out that there are situations in which the notion 
of a strong attractor is not the best choice. 
The notion of a weak attractor has first been considered by 
Ochs~\cite{ochs}.  \par\smallskip

\aussage{Definition}  \label{tq1}
Let~$\varphi$ be an RDS on~$E$ over the MDS 
$\bigl(\Omega,\mathscr F,(\vartheta_t)_{t\in\mathbb R},P\bigr)$. 
Let $\mathcal B\subset2^E$ be an arbitrary subset of the 
power set of~$E$. 
A random set $\omega\mapsto A(\omega)$ is called a 
\emph{weak $\mathcal B$-attractor} for~$\varphi$ if~$A$ 
satisfies conditions~(a) (compactness) and~(b) (almost sure 
strict invariance) of Definition~\ref{sq2}.\ref{f_ibf_watr}, 
and if 
\begin{displaymath}
  \leqno\mbox{(\=c)}\quad\displaystyle
  \lim_{t\to\infty}\sup_{x\in B}
  d\bigl(\varphi(t,\omega)(x),A(\vartheta_t\omega)\bigr)=0 
  \qquad\mbox{in probability for every}\ B\in\mathcal B. 
\end{displaymath}

\aussage{Remark}\label{tq7}
(i) \ The notion of a strong attractor is also often 
referred to as a \emph{pullback attractor}.  \\[.7ex]
(ii) \ Asking for convergence almost surely in (\=c) gives yet 
another concept, referred to as a \emph{forward attractor}. 
Clearly both pullback attractors and forward attractors 
are weak attractors. 
However, a weak attractor need neither be a strong (pullback) 
nor a forward attractor, see Scheutzow~\cite{scheutzow02}. 
Also uniform exponential attraction does not suffice to imply 
that a forward attractor is also a pullback 
attractor, see Crauel~\cite{crauel02a}.  \par\medskip

For further literature dealing with the notion of weak 
attractors 
see Ashwin and Ochs~\cite{ashwin-ochs}, 
Arnold and Schmalfu\ss~\cite{arnold-schmalfuss}, 
Chueshov and Scheutzow~\cite{chueshov-scheutzow}, 
Crauel, Duc and Siegmund~\cite{crauel-duc-siegmund}, 
and Scheutzow~\cite{scheutzow08}.  \par

We will need some further notions. 
For non-empty sets $A,B\subset E$, $E$ a metric space, with~$B$ 
bounded we denote by 
\begin{equation}\label{tq4}
  d(B,A)=\sup_{b\in B}d(b,A)=\sup_{b\in B}\inf_{a\in A}d(a,b)
\end{equation}
the Hausdorff semi-distance. 
It should not cause confusion to use the same letter~$d$ for the 
metric on~$E$ and for the Hausdorff semi-distance on subsets of~$E$. 
The Hausdorff metric between two compact sets $A,B\subset E$ is 
given by $d_{\rm H}(A,B)=\max\{d(A,B),d(B,A)\}$. 
\par\smallskip

\aussage{Definition}  \label{tq2}
Suppose that~$\varphi$ is an RDS on~$E$ over 
$\bigl(\Omega,\mathscr F,(\vartheta_t)_{t\in\mathbb R},P\bigr)$. 
\begin{enumerate}[(i)]
\item A random set $\omega\mapsto K(\omega)$ is said to 
  \emph{attract} another random set $\omega\mapsto B(\omega)$ 
  \emph{strongly}, if 
  \begin{displaymath}
    \lim_{t\to\infty}
    d\bigl(
      \varphi(t,\vartheta_{-t}\omega)B(\vartheta_{-t}\omega), 
      K(\omega)\bigr)=0
      \qquad\mbox{for}\ P\mbox{-almost every}\ \omega\in\Omega. 
  \end{displaymath}
\item $K$ is said to \emph{attract}~$B$ \emph{weakly} if 
  \begin{displaymath}
    \lim_{t\to\infty}
    d\bigl(
      \varphi(t,\vartheta_{-t}\omega)B(\vartheta_{-t}\omega), 
      K(\omega)\bigr)=0
      \qquad\mbox{in probability.} 
  \end{displaymath}
\item The \emph{$\Omega$-limit set} of a random set 
  $\omega\mapsto B(\omega)$ is the random set given by 
  \begin{displaymath}
    \Omega_B(\omega)=\bigcap_{T\geq0}
      \overline{\bigcup_{t\geq T}
      \varphi(t,\vartheta_{-t}\omega)B(\vartheta_{-t}\omega)}. 
  \end{displaymath}
\end{enumerate}

We will make use of invariance of $\Omega$-limit sets. 
A random set $\omega\mapsto B(\omega)$ is said to be 
\emph{$\varphi$-invariant} if for every $t\geq0$ there exists 
a set $\Omega_t$ of full measure, such that 
$\varphi(t,\omega)(B(\omega))\subset B(\vartheta_t\omega)$ 
for all $\omega\in\Omega_t$. 
Compare with the notion of strict invariance introduced in 
Definition~\ref{sq2}.\ref{f_ibf_watr}\,(b). 
Note that these notions are not used consistently in the 
literature. 
Often ``forward invariant/invariant'' is used instead of 
``invariant/strictly invariant''.

\aussage{Lemma}\label{tq3}%
\emph{Any $\Omega$-limit set $\Omega_B$ is $\varphi$-invariant.} 

See Crauel~\cite{crauel99} Lemma~5.1, for the proof.

\section{Criteria for strong attractors}\label{sq3}
In this section, we assume that~$\varphi$ is a continuous RDS 
over the metric dynamical system 
$\big(\Omega,\mathscr F,(\vartheta_t)_{t\in\mathbb R},P\big)$, 
taking values in the Polish space~$E$, which is equipped with 
a metric~$d$. 
We note that we assume continuous time here, mainly to ease 
notation. 
The results hold for discrete time as well, and the proofs 
remain the same. 
It should be emphasized that we do not assume continuity in 
the dependence on time.  \\
For a subset~$A$ of~$E$ we denote the closed $\delta$-neighbourhood 
of~$A$ by $A^\delta$.  \par\smallskip

\aussage{Theorem}\label{strongB}
\emph{The following are equivalent: 
\begin{itemize}
\item[\rm(i)]$\varphi$ has a strong $B$-attractor. 
\item[\rm(ii)]For every $\eps>0$ there exists a 
  compact subset $C_\eps$ such that for each $\delta>0$ 
  and each bounded and closed subset~$B$ of~$E$ it holds that 
  \begin{align*}
    P\bigl\{B\subset\bigcup_{s\geq 0}\bigcap_{t\geq s}
    \varphi(t,\vartheta_{-t}\omega)^{-1}(C_\eps^\delta)\bigr\}
    \geq1-\eps. 
  \end{align*}
\item[\rm(iii)] There exists a compact strongly 
  $B$-attracting set $\omega\mapsto K(\omega)$. 
\end{itemize}}

{\sc Proof} \ 
Equivalence of~(i) and~(iii) is proved in Crauel~\cite{crauel01}, 
Theorem~3.4 and Remark~3.5.  \par

To see that~(i) implies~(ii), fix $\eps>0$. 
Since~$E$ is a Polish space and the attractor~$A$ is a random 
variable taking values in the compact sets, there exists a 
compact subset $C_\eps\subset E$ such that 
\begin{displaymath}
  P\bigl\{A(\omega)\subset C_\eps\bigr\}\geq 1-\eps 
\end{displaymath}
(see Crauel~\cite{crauel02b} Proposition~2.15). 
For $\delta>0$ and a bounded and closed subset~$B$ of~$E$ we have 
\begin{displaymath}
  P\Bigl\{B\subset
  \bigcup_{s\geq 0}\bigcap_{t\geq s}
  \varphi(t,\vartheta_{-t}\omega)^{-1}(A^\delta(\omega))\Bigr\} =1. 
\end{displaymath}
Consequently, 
\begin{eqnarray*}
  \lefteqn{P\Bigl\{
    B\subset
    \bigcup_{s\geq 0}\bigcap_{t\geq s}
    \varphi(t,\vartheta_{-t}\omega)^{-1}(C_\eps^\delta)\Bigr\}}  \\ 
  &\geq&
  P\Bigl\{
  B\subset
  \bigcup_{s\geq 0}\bigcap_{t\geq s}
  \varphi(t,\vartheta_{-t}\omega)^{-1}(A^\delta(\omega))\Bigr\}
  - P\bigl\{A(\omega)\not\subset C_\eps\bigr\}  \\
  &\geq&1-\eps,
\end{eqnarray*}
proving (ii).  \par

Finally, in order to show that~(ii) implies~(i), let 
$B_0\subset B_1\subset\ldots$ be a sequence of bounded and closed 
subsets on~$E$ such that for each bounded set~$B$ there exists 
some~$k\in\mathbb N$ with~$B\subset B_k$. 
For instance, $B_k$ may be taken to be the ball of radius~$k$ 
around some $x_0\in E$. 
Put 
\begin{displaymath}
  A(\omega)=\overline{\bigcup_{k\in\mathbb N}\Omega_{B_k}(\omega)},  
\end{displaymath}
where $\Omega_{B_k}(\omega)$ denotes the $\Omega$-limit set of~$B_k$ 
(see Definition~\ref{sq2}.\ref{tq2}\,(iii)). 
Let $\eps>0$. 
By condition~(ii) there exists a compact set~$C_\eps\subset E$ such 
that for all~$k$ and all $\delta >0$ we have 
\begin{displaymath}
  P\bigl\{\Omega_{B_k}(\omega)\subset C_\eps^\delta\bigr\}
  \geq 1-\eps, 
\end{displaymath}
which implies $P\{A(\omega)\subset C_\eps\}\geq 1-\eps$. 
In particular, $A(\omega)$ as well as $\Omega_{B_k}(\omega)$ is 
compact for every~$k$, for $P$-almost every $\omega\in\Omega$.  \par

We now claim that $\Omega_{B_k}$ is strictly invariant for 
every~$k$. 
Invariance of~$\Omega_{B_k}$ follows from Lemma~\ref{sq2}.\ref{tq3}. 
In order to see that $\Omega_{B_k}$ is strictly invariant, fix $\eps>0$ 
and $t\geq0$, and suppose that $y\in\Omega_{B_k}(\vartheta_t\omega)$. 
Then $y=\lim\varphi\bigl(t_n,\vartheta_{-t_n}(\vartheta_t\omega)\bigr)b_n$ 
for some sequence $t_n\to\infty$, $b_n\in B_k$, with $n\to\infty$. 
Consider the sequence $\varphi(t_n-t,\vartheta_{-(t-t_n)}\omega)b_n$, 
defined for~$n$ with $t_n-t\geq0$. 
Then for~$\omega$ with $\Omega_{B_k}(\omega)\subset C_\eps$ we have 
$\lim d\bigl(\varphi(t_n-t,\vartheta_{-(t_n-t)}\omega)b_n,C_\eps\bigr)=0$ 
for $n\to\infty$. 
This implies existence of a convergent subsequence with limit~$z(\omega)$, 
say. 
Clearly $z(\omega)\in\Omega_{B_k}$. 
Using the same notation for the subsequence, continuity of~$\varphi(t,\omega)$ 
implies
\begin{displaymath}
  \varphi(t,\omega)z(\omega)
  =\lim_{n\to\infty}\varphi(t_n,\vartheta_{-(t_n-t)}\omega)b_n=y(\omega). 
\end{displaymath}
Consequently, 
\begin{displaymath}
  P\bigl\{\Omega_{B_k}(\vartheta_t\omega)
  \subset\varphi(t,\omega)\Omega_{B_k}(\omega)\bigr\}\geq1-\eps. 
\end{displaymath}
This holding for every $\eps>0$, we obtain that $\Omega_{B_k}$ is 
strictly invariant, whence also~$A$ is strictly invariant.  \par

It remains to show that~$A$ attracts every bounded set. 
It suffices to show that~$A$ attracts each of the sets $B_k$. 
Now $\Omega_{B_k}$ attracts~$B_k$ for every~$k$, $P$-almost surely. 
In fact, if this would not be the case then there would be a 
$\delta>0$, a sequence $t_n\to\infty$, and $b_n\in B_k$ such that 
$d\bigl(\varphi(t_n,\vartheta_{-t_n}\omega)b_n,\Omega_{B_k}(\omega)\bigr)
\geq\delta$ for every $n\in\mathbb N$ with positive probability. 
In view of the fact that 
$\bigl(\varphi(t_n,\vartheta_{-t_n}\omega)b_n\bigr)_{n\in\mathbb N}$ 
has a convergent subsequence with probability larger than $1-\eps$ 
for $\eps>0$ arbitrary, which then must converge to a limit in 
$\Omega_{B_k}(\omega)$, this would yield a contradiction, completing 
the proof.  \qed

Next we consider the case of a strong $C$-attractor. 
Several arguments of the previous Theorem~\ref{sq3}.\ref{strongB} 
proceed analogously. 
However, the construction of the attractor has to be done in a 
different way, and also the argument invoked for the verification 
of the attraction property is different. 
Therefore it seems appropriate to formulate the result separately, 
even though the assertions appear to be very similar.  \par\smallskip

\aussage{Theorem}\label{strongC}
\emph{The following are equivalent: 
  \begin{itemize}
  \item[\rm(i)]$\varphi$ has a strong $C$-attractor.
  \item[\rm(ii)]For every $\eps>0$ there exists a 
    compact subset $C_\eps$ such that for each $\delta>0$ 
    and each compact subset~$K$ of~$E$ it holds that 
    \begin{displaymath}
      P\bigl\{K\subset
      \bigcup_{s\geq 0}
      \bigcap_{t\geq s}
      \varphi(t,\vartheta_{-t}\omega)^{-1}(C_\eps^\delta)\bigr\}
      \geq 1-\eps.
    \end{displaymath}
  \item[\rm(iii)] There exists a compact strongly 
    $C$-attracting set $\omega\mapsto K(\omega)$. 
  \end{itemize}}

{\sc Proof} \ 
Equivalence of~(i) and~(iii) follows again from Crauel~\cite{crauel01}, 
Theorem~3.4 and Remark~3.5, and (i)$\Rightarrow$(ii) is proved 
exactly as in Theorem~\ref{sq3}.\ref{strongB}.  \par

In order to see that~(ii) implies~(i), 
consider $C_{1/k}$, $k\in\mathbb N$, where the sets $C_\eps$ 
are from condition~(ii), and put 
\begin{displaymath}
  A(\omega)=\overline{\bigcup_{k\in\mathbb N}\Omega_{C_{1/k}}(\omega)}. 
\end{displaymath}
As in Theorem~\ref{sq3}.\ref{strongB} we obtain that~$A$ is 
compact almost surely, and strictly invariant. 
It remains to show that~$A$ attracts every compact set. 
Let~$K\subset E$ be compact. 
Then $\Omega_K$ is strictly invariant, arguing as in the proof 
of Theorem~\ref{sq3}.\ref{strongB}. 
Further 
\begin{displaymath}
  1-\frac1n\leq
  P\bigl\{\Omega_K(\omega)\subset C_{1/n}\bigr\}
  \leq P\bigl\{\Omega_K(\omega) \subset \Omega_{C_{1/n}}\bigr\} 
  \qquad\mbox{for every}\ n\in\mathbb N 
\end{displaymath}
(see Crauel~\cite{crauel99} Proposition~5.2 for the second 
inequality), hence $\Omega_K\subset A$ almost surely, and 
therefore~$A$ attracts~$K$ almost surely. 
This completes the proof.  \qed

\section{Criteria for weak attractors}\label{weak}

In this section we are interested in weak attractors. 
We follow the same line as in the previous section on strong attractors. 
We establish necessary and sufficient criteria first for the existence 
of weak $B$-attractors, and then for the existence of weak 
$C$-attractors. 
The structure is very similar to that of Theorems~\ref{sq3}.\ref{strongB} 
and~\ref{sq3}.\ref{strongC}. 
The difference of the two concepts gets visible in the corresponding 
conditions~(ii), which is eventually ``uniform in time'' for strong 
attractors, while it is eventually ``pointwise in time'' for weak 
attractors.  \par

Again we consider the cases of $B$-attractors and of $C$-attractors 
separately, even if the assertions are very similar. 
Also certain parts of the proofs are similar, and they could be 
presented in one result. 
However, the construction of the attractor is different in those 
two cases. 
The argument is more straightforward for $B$-attractors. 
Therefore the presentation has been split, and the arguments are 
given separately.  \par

As before,~$\varphi$ is a continuous RDS over the metric dynamical system 
$\big(\Omega,\mathscr F,(\vartheta_t)_{t\in\mathbb R},P\big)$, taking 
values in the Polish space $E$, equipped with a metric~$d$.  \par

We will make use of an elementary lemma. 
Again~$B^\delta$ denotes the $\delta$-neighbourhood of a subset~$B$ of 
a metric space.  \par\smallskip

\aussage{Lemma}\label{tq5}
\emph{Suppose that $\varphi:X\to Y$ is a continuous map from a 
  metric space~$X$ to a metric space~$Y$. 
  Let $C\subset X$ be compact.}
\begin{enumerate}[(i)]
\item\emph{For every $\eps>0$ there exists a $\gamma>0$ such that 
    \begin{displaymath}
      \varphi(C^\gamma)\subset\bigl(\varphi(C)\bigr)^\eps, 
    \end{displaymath}
    or, equivalently, $d\bigl(\varphi(C^\delta),\varphi(C)\bigr)\to0$ 
    for $\delta\to0$.}
\item\emph{If $B\subset Y$ satisfies $B\subset\overline{\varphi(C^\delta)}$ 
    for every $\delta>0$, then $B\subset\varphi(C)$.}
\end{enumerate}
  
{\sc Proof} \ Assuming~(i) not to be true, we get existence of some $\eps>0$ 
and of a sequence $(x_n)_{n\in\mathbb N}$ with $x_n\in C^{1/n}$ such that 
\begin{equation}\label{tq6}
  d\bigl(\varphi(x_n),\varphi(C)\bigr)\geq\eps
  \qquad\mbox{for every}\ n\in\mathbb N. 
\end{equation}
By compactness of~$C$ the sequence $(x_n)$ has a convergent 
subsequence, denoted by $(x_n)$ again, converging to some~$x\in C$. 
This implies convergence of $d\bigl(\varphi(x_n),\varphi(C)\bigr)$ to 
$d\bigl(\varphi(x),\varphi(C)\bigr)=0$, contradicting~\eqref{tq6}.  \par

In order to obtain~(ii) note that $B\subset\overline{\varphi(C^\delta)}$ implies 
$d\bigl(B,\varphi(C)\bigr)
  \leq d\bigl(\overline{\varphi(C^\delta)},\varphi(C)\bigr)
  =d\bigl(\varphi(C^\delta),\varphi(C)\bigr)$, 
which converges to zero for~$\delta\to0$ by~(i). 
Therefore $d\bigl(B,\varphi(C)\bigr)=0$, which implies $B\subset\varphi(C)$ 
since~$\varphi(C)$ is closed, in fact even compact.  \qed


\aussage{Theorem}\label{weakB}
\emph{The following are equivalent: 
  \begin{itemize}
  \item[\rm(i)] $\varphi$ has a weak $B$-attractor. 
  \item[\rm(ii)] For every  $\eps>0$ there exists a compact 
    subset $C_\eps$ such that for each $\delta > 0$ and each 
    bounded subset $B$ of $E$ there is a $t_0>0$ with the 
    property that for all $t\geq t_0$, 
    \begin{displaymath}
      P\bigl\{\varphi(t,\omega)(B)\subset C_\eps^\delta\bigr\}\geq 1-\eps.
    \end{displaymath}
  \item[\rm(iii)] There exists a compact weakly $B$-attracting 
    set $\omega\mapsto K(\omega)$. 
  \end{itemize}}

{\sc Proof} \ 
Obviously~(i) implies~(iii).  \par

To see that~(iii) implies~(ii), let $\omega\mapsto K(\omega)$ be 
a compact weakly $B$-attracting set. 
Fix $\eps>0$. 
Since~$E$ is a Polish space and~$K$ is a random variable taking values 
in the compact sets, there exists a compact subset~$C_\eps$ of~$E$ 
such that 
\begin{displaymath}
  P\bigl\{K(\omega)\subset C_\eps\bigr\}\geq1-\frac\eps2
\end{displaymath}
(see Crauel~\cite{crauel02b} Proposition~2.15). 
For every $\delta>0$ and every bounded subset~$B$ of~$E$ there 
exists~$t_0$ such that for every $t\geq t_0$ one has 
\begin{displaymath}
  P\bigl\{\varphi(t,\vartheta_{-t}\omega)(B)\subset K^\delta(\omega)\bigr\}
  \geq1-\frac\eps2. 
\end{displaymath}
Therefore, for every $t\geq t_0$ 
\begin{align*}
  P\bigl\{\varphi(t,\omega)(B)\not\subset C_\eps^\delta\bigr\}
  &=
  P\bigl\{
  \varphi(t,\vartheta_{-t}\omega)(B)\not\subset C_\eps^\delta\bigr\}  \\
  &\leq
  P\bigl\{\varphi(t,\vartheta_{-t}\omega)(B)\not\subset K^\delta(\omega)\bigr\}
  +P\bigl\{K(\omega)\not\subset C_\eps\bigr\}  \\
  &\leq\eps,
\end{align*}
which proves~(ii).  \par

Finally, let us show that~(ii) implies~(i). 
Fix a point $x_0\in E$ and let $B_k$, $k\in\mathbb N$, be closed balls 
in~$E$ with center $x_0$ and radii increasing to infinity, such that 
$C_{2^{-k}}^1\subset B_k$, $k\in\mathbb N$. 
Define a sequence of numbers $u_n>n$, $n\in\mathbb N$, recursively by 
\begin{equation}\label{eins}
  P\bigl\{\varphi({\textstyle\sum}_{k=1}^n u_k,\omega)(B_n) 
  \subset C_{2^{-m}}^{1/n}\bigr\} \geq 1-2^{-m},\;m=1,\ldots,n, 
\end{equation}
and
\begin{equation}\label{zwei}
  P\{ \varphi(u-n,\omega)(B_n) \subset B_{n-1}\} 
  \geq 1 - 2^{-n+1}\ \mbox{for all}\ u\geq u_n,\;n \geq 2. 
\end{equation}
Define $t_n=\sum_{i=1}^n u_i$, $n\in\mathbb N$, and put 
\begin{displaymath}
  A(\omega)=
  \bigcap_{n\in\mathbb N}
  \overline{\bigcup_{k\geq n}\varphi(t_k,\vartheta_{-t_k}\omega)(B_k)}.
\end{displaymath}
Then~\eqref{zwei} implies 
\begin{displaymath}
  \sum_{k=2}^\infty
    P\bigl\{\varphi(u_k,\omega)(B_k)\not\subset B_{k-1}\bigr\}<\infty, 
\end{displaymath}
so that the first Borel-Cantelli lemma yields existence of a 
set~$\Omega_0$ of full measure and a positive integer $j_0(\omega)$ 
such that 
\begin{equation}\label{BC}
  \varphi\bigl(u_j,\vartheta_{-t_j}\omega\bigr)(B_j)\subset B_{j-1}
\end{equation}
for every $j\geq j_0(\omega)$, $\omega\in\Omega_0$. 
Here we also have made use of the $\vartheta_t$-invariance of~$P$ 
for every fixed $t\in\mathbb R$. 
In particular, using the cocycle property, we have 
\begin{displaymath}
  A(\omega)=
  \bigcap_{j=j_0(\omega)}^\infty
  \overline{\varphi(t_j,\vartheta_{-t_j}\omega)(B_j)}\qquad 
  \mbox{on}\ \Omega_0.
\end{displaymath}

We claim that~$A$ is a weak $B$-attractor.  \par

{\sc Step~1} \ $A$ is almost surely compact.  \\
Fix $m\in\mathbb N$. 
Using~\eqref{eins} we get, for $n\geq m$, 
\begin{eqnarray*}
  \lefteqn{P\bigl\{A(\omega)\subset C_{2^{-m}}^{1/n}\bigr\}}  \\
  &\geq&P\bigl\{
  \overline{\varphi(t_n,\vartheta_{-t_n}\omega)(B_n)}\subset C_{2^{-m}}^{1/n}\}
  -P\{j_0(\omega)>n\}  \\
  &\geq&1-2^{-m}-P\{j_0(\omega)>n\}\to1-2^{-m}\ \mbox{for}\ n\to\infty, 
\end{eqnarray*}
hence $P\{A(\omega)\subset C_{2^{-m}}\}\geq1-2^{-m}$ for every $m\in\mathbb N$, 
so that~$A$ is almost surely compact.  \par

{\sc Step~2} \ $A$~is strictly invariant.  \\
Fix $t>0$. 
We first establish $\varphi(t,\omega)A(\omega)\subset A(\vartheta_t\omega)$ 
almost surely.  \\
Fix $\eps>0$ and $\delta>0$ and choose~$n$ so large that 
\begin{enumerate}[(i)]
\item $\displaystyle P\{\varphi(t_n,\vartheta_{-t_n}\omega)(B_n) 
  \subset A^\delta(\omega)\}\geq1-\frac\eps3$,
\item $\displaystyle P\{ j_0(\omega)>n+1\}\leq\frac\eps3$, 
\item $\displaystyle 2^{-n}\leq\frac\eps3$. 
\end{enumerate}
Observe that such an~$n$ always exists. 
Then 
\begin{eqnarray*}
  \lefteqn{P\{\varphi(t,\omega)A(\omega)
    \subset A^\delta(\vartheta_t\omega)\}}  \\
  &\geq&
  P\bigl\{
  \varphi(t,\omega)\varphi(t_{n+1},\vartheta_{-t_{n+1}}\omega)(B_{n+1}) 
  \subset A^\delta(\vartheta_t\omega)\bigr\}
  -P\{j_0(\omega)>n+1\}  \\
  &\geq&
  P\bigl\{\varphi(t+t_{n+1},\vartheta_{-t_{n+1}}\omega) (B_{n+1})
  \subset\varphi(t_n,\vartheta_{-t_n}\vartheta_t\omega)(B_n)\bigr\}
  -\frac\eps3-P\{j_0(\omega)>n+1\}  \\
  &\geq&
  P\bigl\{\varphi(t+t_{n+1}-t_n,\omega)(B_{n+1})\subset B_n\bigr\}
  -\frac{2\eps}3  \\
  &\geq&1-\eps, 
\end{eqnarray*}
where we used~(i) and the cocycle property for the second inequality, 
(ii) and the cocycle property for the third inequality, 
and~\eqref{zwei} and~(iii) for the last inequality. 
Since $\delta,\eps>0$ are arbitrary, this implies 
$\varphi(t,\omega)A(\omega)\subset A(\vartheta_t\omega)$ 
almost surely.  \par

Next we prove 
$A(\vartheta_t\omega)\subset\varphi(t,\omega)A(\omega)$ almost surely. 
Fix $\eps>0$ and $\delta>0$ and choose~$i$ so large that 
\begin{enumerate}[(i)]
\item $\displaystyle P\bigl\{
  \varphi(t_i,\vartheta_{-t_i}\omega)(B_i)
  \subset A^\delta(\omega)\bigr\}\geq 1-\frac\eps3$, 
  and, consequently, 
  \begin{equation}\label{tq8}
    P\bigl\{\overline{\varphi(t,\omega)\varphi(t_i,\vartheta_{-t_i}\omega)(B_i)}
    \subset\overline{\varphi(t,\omega)A^\delta(\omega)}\bigr\}
  \geq 1-\frac\eps3
  \end{equation}
\item $\displaystyle P\{j_0(\omega)>i+1\}\leq\frac\eps3$
\item $\displaystyle2^{-i}\leq\frac\eps3$
\item $i \geq t-1$.
\end{enumerate}
Observe that such an~$i$ always exists. 
Then 
\begin{eqnarray*}
  \lefteqn{P\bigl\{
    A(\vartheta_t\omega)
      \subset\overline{\varphi(t,\omega)A^\delta(\omega)}\bigr\}}  \\
  &\geq&P\Bigl\{
    A(\vartheta_t\omega)
      \subset
      \overline{
        \varphi(t,\omega)\varphi(t_i,\vartheta_{-t_i}\omega)(B_i)}\Bigr\}
    -\frac\eps3  \\
  &\geq&
  P\Bigl\{\overline{\varphi(t_{i+1},\vartheta_{t-t_{i+1}}\omega)(B_{i+1})}
    \subset
    \overline{
      \varphi(t,\omega)\varphi(t_i,\vartheta_{-t_i}\omega)(B_i)}\Bigr\}
  -P\{j_0(\omega)>i+1\}-\frac\eps3  \\
  &\geq&
  P\bigl\{\varphi(t_{i+1},\vartheta_{t-t_{i+1}}\omega)(B_{i+1})
  \subset\varphi(t,\omega)\varphi(t_i,\vartheta_{-t_i}\omega)(B_i)\bigr\}
  -P\{j_0(\omega)>i+1\}-\frac\eps3  \\
  &\geq&
  P\bigl\{\varphi(t_{i+1}-t_i-t,\vartheta_{t-t_{i+1}}\omega)(B_{i+1})
  \subset B_i\bigr\}
  -P\{j_0(\omega)>i+1\}-\frac\eps3  \\
  &\geq&1-\eps,
\end{eqnarray*}
where we used~\eqref{tq8} for the first inequality, the cocycle property 
for the fourth inequality, and~(ii), (iii), (iv) as well as~\eqref{zwei} 
for the final inequality. 
Since $\eps>0$ is arbitrary, we get 
$A(\vartheta_t\omega)\subset\overline{\varphi(t,\omega)A^\delta(\omega)}$ 
almost surely for every $\delta>0$, which, by virtue of 
Lemma~\ref{weak}.\ref{tq5}\,(ii), implies 
$A(\vartheta_t\omega)\subset\varphi(t,\omega)A(\omega)$ almost surely.  \par

{\sc Step~3} \ $A$ attracts every bounded $B\subset E$ in 
probability.  \\
Let~$B$ be a bounded subset of~$E$, and fix $\delta,\eps>0$. 
Then there exists~$j$ such that $2^{-j}\leq\eps/2$, 
$B\subset B_{j+1}$, and 
$P\{\varphi(t_j,\vartheta_{-t_j}\omega)(B_j)
  \subset A^\delta(\omega)\}\geq 1-\eps/2$. 
Then, for every $t\geq t_{j+1}$, 
\begin{eqnarray*}
  \lefteqn{P\bigl\{
    \varphi(t,\vartheta_{-t}\omega)(B)\subset A^\delta(\omega)\bigr\}}  \\
  &\geq&
  P\bigl\{
  \varphi(t_j,\vartheta_{-t_j}\omega)(B_j)\subset A^\delta(\omega)\bigr\}
  -P\bigl\{
  \varphi(t,\vartheta_{-t}\omega)(B_{j+1})
    \not\subset\varphi(t_j,\vartheta_{-t_j}\omega)(B_j)\bigr\}  \\
  &\geq&1-\eps,
\end{eqnarray*}
where we used~\eqref{zwei} for the final inequality.  \qed

In order to obtain a corresponding result for weak $C$-attractors we 
will make use of a technical lemma.  \par

\aussage{Lemma}\label{vvv}
\emph{Let $\varphi:[0,\infty)\times\Omega\times E\to E$ be as before. 
Let~$C\subset E$ be compact and~$B\subset E$ closed. 
Let $\alpha>0$. 
Then the map 
\begin{displaymath}
  t\mapsto\gamma(t)
  :=\sup\bigl\{\eta>0: P\{\varphi(t,\omega)C^\eta\subset B\}\geq\alpha\bigr\}
\end{displaymath}
is measurable.}

{\sc Proof} \
First note that for every $r\in\mathbb[0,\infty)$ the set 
$\{(t,\omega)\in[0,\infty)\times\Omega:\varphi(t,\omega)C^r\subset B\}$ 
is measurable due to measurability of $(t,\omega)\mapsto\varphi(t,\omega)x$ 
for every $x\in E$ together with separability of the metrizable~$E$, which 
implies separability of every subset of~$E$. 
Put 
\begin{displaymath}
  V=\bigcup_{q\in\mathbb Q^+}
    \Bigl(\bigl\{(t,\omega):\varphi(t,\omega)C^q\subset B\bigr\}\times[0,q)
    \Bigr), 
\end{displaymath}
then $V$ is a measurable subset of $[0,\infty)\times\Omega\times[0,\infty)$. 
Therefore $(t,\eta)\mapsto f(t,\eta)=\int 1_V(t,\omega,\eta)\,\mathrm dP(\omega)$ 
is measurable, and so is 
$W=\{(t,\eta):f(t,\eta)\geq\alpha\}\subset[0,\infty)\times[0,\infty)$ 
and, consequently, also 
$t\mapsto\gamma(t)=\int1_W(t,\eta)\,\mathrm d\eta$.  \qed

\aussage{Theorem}\label{weakC}
\emph{The following are equivalent:
\begin{itemize}
\item[\rm(i)] $\varphi$ has a weak $C$-attractor.
\item[\rm(ii)] For every $\eps>0$ there exists a compact subset~$C_\eps$
  such that for each $\delta>0$ and each compact subset~$C$ of~$E$ there
  is a $t_0>0$ with the property that for all $t\geq t_0$
  \begin{displaymath}
    P\bigl\{\varphi(t,\omega)(C)\subset C_\eps^\delta\bigr\}\geq 1-\eps.
  \end{displaymath}
\item[\rm(iii)] There exists a compact weakly $C$-attracting set
  $\omega\mapsto K(\omega)$.
\end{itemize}}

{\sc Proof} \
Obviously~(i) implies~(iii), and~(iii) implies~(ii) by exactly
the same argument as in the proof of Theorem~\ref{weak}.\ref{weakB}.  \par

In order to prove that~(ii) implies~(i) we will follow the proof of 
Theorem~\ref{weak}.\ref{weakB} as closely as possible. 
Lemma~\ref{weak}.\ref{tq5} and Lemma~\ref{weak}.\ref{vvv} imply that for 
every strictly positive sequence~$\delta_n$, $n \in \mathbb N$, there exist 
strictly positive sequences~$\gamma_n$, and $u_n>n$ together with 
measurable subsets~$U_n$ of $[u_n/2,2 u_n]$, $n\in\mathbb N$, such that 
$u_n\in U_n$ and~$U_n$ has Lebesgue measure at least $\frac32u_n-\delta_n$ 
for all $n\in\mathbb N$, and such that the sets 
$B_n:=C_{2^{-n}}^{\gamma_n}$, $n \in \mathbb N$, satisfy 
\begin{equation}\label{einsC}
  P\bigl\{\varphi({\textstyle\sum}_{i=1}^n u_i,\omega)(B_n)
  \subset C_{2^{-m}}^{1/n}\bigr\} \geq 1-2^{-m+1},\;m=1,\ldots,n,
\end{equation}
and
\begin{equation}\label{zweiC}
  P\bigl\{\varphi(u,\omega)(B_n) \subset B_{n-1}\bigr\}
  \geq 1-2^{-n+2}\ 
  \mbox{for all}\ u\in U_n,\;n\geq2. 
\end{equation}
Now choose a summable sequence $(\delta_n)$, 
define $t_n=\sum_{i=1}^n u_i$, $n\in\mathbb N$, and put 
\begin{displaymath}
  A(\omega)=
  \bigcap_{n\in\mathbb N}
  \overline{\bigcup_{k\geq n}\varphi(t_k,\vartheta_{-t_k}\omega)(B_k)}.
\end{displaymath}
Defining~$\Omega_0$ and~$\omega\mapsto j_0(\omega)$ as in the proof of 
Theorem~\ref{weak}.\ref{weakB} it follows that 
\begin{equation}  \label{c-attractor1} 
  A(\omega)=
  \bigcap_{j=j_0(\omega)}^\infty
  \overline{\varphi(t_j,\vartheta_{-t_j}\omega)(B_j)}\qquad
  \mbox{on}\ \Omega_0.
\end{equation}

We claim that~$A$ is a weak $C$-attractor.  \par

Almost sure compactness of $A$ follows like in the proof of 
Theorem~\ref{weak}.\ref{weakB}.  \par

To show strict invariance we first replace both conditions~(iii) in 
the proof of Step 2 of Theorem~\ref{weak}.\ref{weakB} by multiplying 
the left hand side by~$2$ (this is due to the factor~$2$ in~\eqref{zweiC} 
compared to~\eqref{zwei}). 
Then all estimates follow as before except that we require that 
$t+u_{n+1}\in U_{n+1}$ respectively $u_{i+1}-t\in U_{i+1}$. 
Due to the definition of the sets~$U_n$ and the summability of 
the sequence~$(\delta_n)$ it follows that for Lebesgue almost 
all $t>0$ we have $t+u_n\in U_n$ and $u_n-t\in U_n$ for all but finitely 
many $n$. 
Therefore we see that the set 
\begin{displaymath}
  T:=\{t \ge 0: 
  \varphi(t,\omega) A (\omega) = A (\vartheta_t \omega)\ \mbox{almost surely}\}
\end{displaymath}
has full Lebesgue measure (and contains 0). 
Further the cocycle property shows that~$T$ is closed under addition, 
from which we conclude that $T=[0,\infty)$.  \par

It remains to prove that~$A$ attracts compact sets.  \\
For $\delta,\,\eps>0$ choose $k_0\in\mathbb N$ with $2^{k_0}>2/\eps$ 
such that for all $k\geq k_0$ 
\begin{equation}\label{tq9}
  P\bigl(
  d(\varphi(t_k,\vartheta_{-t_k}\omega)( B_k),A(\omega))>\delta\bigr)
  <\frac\eps2,
\end{equation}
this is possible due to~\eqref{c-attractor1}. 
Let $C\subset E$ be compact. 
Since~$B_k=C_{2^{-k}}^{\gamma_k}$, condition~(ii) yields existence of 
a time~$t_C>0$ such that for all $t\geq t_C$ 
\begin{displaymath}
  P\bigl\{\varphi(t,\omega)(C)\subset B_k\bigr\}\geq 1-\frac\eps2.
\end{displaymath}
Using the cocycle property we get for all $t\geq t_C$ 
\begin{eqnarray*}
  1-\frac\eps2
  &\leq&
  P\bigl\{
  \varphi(t,\vartheta_{-t_k-t}\omega)(C)\subset B_k\bigr\}  \\
  &\leq&
  P\bigl\{
  \varphi(t_k,\vartheta_{-t_k}\omega)\varphi(t,\vartheta_{-t_k-t}\omega)(C)
    \subset\varphi(t_k,\vartheta_{-t_k}\omega)( B_k)\bigr\}  \\
  &=&
  P\bigl\{\varphi(t+t_k,\vartheta_{-t_k-t}\omega)(C)
    \subset\varphi(t_k,\vartheta_{-t_k}\omega)(B_k)\bigr\}. 
\end{eqnarray*}
Together with~\eqref{tq9} this implies, for all $t\geq t_C+t_k$, 
\begin{displaymath}
  P\bigl(
    d\bigl(\varphi(t,\vartheta_{-t}\omega)(C),A(\omega)\bigr)>\delta
   \bigr)
  <\eps. 
\end{displaymath}
This holding true for every $\delta>0$ and $\eps>0$ we conclude 
that $d\bigl(\varphi(t,\vartheta_{-t}\omega)(C),A(\omega)\bigr)$ 
converges to zero in probability.   \qed

Finally, we give a necessary and sufficient condition for a 
weak attractor to be also a strong attractor. 
Let $\mathcal B$ be an arbitrary family of deterministic 
subsets of $E$. 

\aussage{Proposition} 
\emph{Suppose that~$A$ is a weak $\mathcal B$-attractor. 
Then the following are equivalent: 
\begin{itemize}
 \item [(i)] $A$ is a strong $\mathcal B$-attractor. 
 \item [(ii)] For every $B\in\mathcal B$ it holds that 
   $\Omega_B\subset A$ almost surely. 
\end{itemize}}

{\sc Proof} \ 
\noindent (i) $\Rightarrow$ (ii): obvious.

\noindent (ii) $\Rightarrow$ (i): let $B\in\mathcal B$. 
Due to~(ii), $A$ attracts~$B$ strongly. 
Since~$A$ is a weak attractor by assumption, (i) follows. 
\qed

\subsection*{Acknowledgement}
Part of this work was done during a stay of HC and MS at the 
Institut Mittag-Leffler, which is gratefully acknowledged.

\end{document}